\documentclass[11pt]{amsart}
\usepackage{amsmath,amssymb}
\usepackage{amsrefs}

\title{Coloring the  600 Cell}
\author{Steve Fisk}
\date{\today}

\newcommand{\tgroup}{$\mathcal{T}$}
\newcommand{\igroup}{$\mathcal{I}$} 
\newcommand{\six}{{\sf\bf S}}	
\newcommand{\three}{{\sf\bf P}}	

\newcommand{\btwosix}{$B^2($\six$)$}

\begin{document}
\begin{abstract}
The 600 cell \six\ has exactly 10 5-colorings.  From these colorings
we can construct the space of colorings $B($\six$)$.  This complex has
1344 colorings, and is isomorphic to the space of 5 by 5 Latin
Squares. These simplices split into 4 copies of a quotient of \six\ by
an involution, and two copies of a space made up of even Latin Squares.
\end{abstract}
\maketitle

\subsection*{Introduction}
The five-colorings of the 600 cell have a surprising structure.  In
order to define this structure, we recall \cite{fisk} the definition
of the space $B(X)$ of five-colorings of a 4-complex $X$:
\begin{quote} 
Consider a five-coloring of $X$ to be a map from $X$ to the colors
$\{1,2,3,4,5\}$. The vertices of $B(X)$ are all the distinct sets
$f^{-1}(v)$, where $v$ is a color, and $f$ is a coloring.  Every
five-coloring $f$ determines a 4-simplex 
$\{f^{-1}(1),\dots,f^{-1}(5)\}$, and all 4-simplices of $B(X)$ are of
this form.
\end{quote}

Since $B(X)$ is a 4-complex, we can once again compute the space $B^2(X)$ of
its five-colorings.  Recall the definition of the map 
$\varphi:X\longrightarrow B^2(X)$:
\begin{quote} 
If $p$ is a vertex of $X$, then $\varphi(p)$ consists of all vertices
$f^{-1}(v)$ such that $p\in f^{-1}(v)$.  
\end{quote}

This is a well-defined simplicial map that sends 4-simplices to
4-simplices.  $\varphi$ is sometimes an isomorphism, but usually is
neither 1-1 nor onto. 

\subsection*{The 600 cell} 

The 600 cell \six\ is the triangulation of the 3-sphere that has 120
vertices, 600 tetrahedra  and the link of every vertex is an icosahedron. 
\six\ can be realized in ${\Bbb R}^4$ in such a way that all tetrahedra
are regular \cite{duval}.

\subsection*{\six\ has 10 five-colorings} 
It is very unusual to be able to determine all the five-colorings of a
graph with 120 vertices.  We are successful  because there is an
inductive way of constructing \six, and  five-colorings extend uniquely
from one stage to the next. 
First, pick any vertex $v$, and fix a five-coloring $f$
of the star $S_0$ of $v$.  $S_0$ is $v$ joined to an icosahedron.  If
$S_1$ is $S_0$ along with all tetrahedra that meet $S_0$, then $f$ has
a unique extension to $S_1$.  Next, if $S_2$ is $S_1$ along with all
tetrahedra meeting $S_1$, then $f$ has a unique extension to $S_2$.
Continuing, $f$ has a unique extension to successive shells, and $f$
extends uniquely to \six.  These assertions about unique extensions
are the result of computer calculations.

  Since the icosahedron has exactly 10 colorings, \six\ does as well.

\subsection*{B(\six) has 25 vertices}

To find the vertices of B(\six) we need an algebraic description of
\six\ and  its colorings.  \six\ can be realized as the convex hull
of the 120 points in ${\Bbb R}^4$ that lie in the quaternionic
icosahedral group \igroup \cite{duval}.  Let \tgroup\ be the subgroup
of \igroup\ whose 24 vertices are the vertices of a 24 cell embedded
in \six.  If $p$ is an element of order 5 in \igroup, then the five
cosets \tgroup$p^k$ (k = 0,1,2,3,4) are all disjoint.  Similarly, the
five cosets $p^k$\tgroup (k=0,1,2,3,4) are all disjoint.  Since each
coset has 24 elements, $5\cdot24 = 120$, and no two of the elements of
\tgroup are adjacent in \six, it follows that we have 2 colorings of
\six, namely $\{$\tgroup$p^k,k=1,\dots,5\}$ and
$\{p^k$\tgroup$,k=1,\dots,5\}$. We can construct 8 additional
colorings by multiplying each of these colorings by $p^1,p^2,p^3$, or
$p^4$ on the appropriate side.

  Since we established computationally that there are exactly 10
colorings, it follows that these are all the colorings.  Furthermore,
we can identify the set of points of one color in these colorings with
pairs $(i,j)$ representing the double coset $p^i$\tgroup$p^j$. A set of
five pairs $(i,j)$ is a coloring if they all have the same first
coordinate or the same second, so B(\six) is isomorphic to $\Delta^4
\# \Delta^4$, the Cartesian product of two 4-simplices
\cite{fisk}. This means that we can represent B(\six) by   the
following diagram, where the vertices are the ``.'''s, and the ``.'''s
in a row (or column) are the 4-simplices.

$$
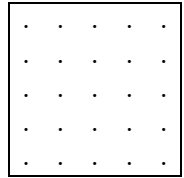
\begin{array}{|ccccc|}
\hline
. & . & . & . & . \\
. & . & . & . & . \\
. & . & . & . & . \\
. & . & . & . & . \\
. & . & . & . & .\\\hline
\end{array}
$$

\subsection*{The 4-simplices of B$^2$(\six) are Latin Squares}
A five coloring of B(\six) is an assignment of one of the integers
$1,2,3,4,5$ to each ``.'' of the figure such that each row and each
column has all distinct entries.  Such labelings are exactly Latin
Squares.  It is easy to compute that there are 1344 5 by 5 Latin
Squares, so B$^2$(\six) has 1344 4-simplices.

\subsection*{The vertices of B$^2$(\six) are permutations}
A  vertex  of  B$^2$(\six) is the set of ``.'' that have the same
label in a Latin Square. Such a set has one ``.'' in each row and in
each column. These sets are  permutations of 
$\{1,2,3,4,5\}$, so there are 120 vertices.

\subsection*{$B^3$(\six) = $B$(\six)} 
In \cite{fisk}, we showed that $S_5$, the space of 5 by 5 Latin
Squares, satisfies $B^2(S_5) = S_5$.  
Since we observed that \btwosix$ = S_5$, the result follows.

\subsection*{The map $\varphi:$\six$\longrightarrow $\btwosix\ is not 1-1.} 
Both \six\ and \tgroup\ are fixed by the involution sending x to -x.
Consequently, whenever a vertex v of \six\ is in some 24-cell, -v is
in the same 24-cell.  It follows that $\varphi(v) = \varphi(-v)$, and
that $\varphi$ is not 1-1 on \six.  Computationally (or by using the
quaternionic representation), we find that -v
is the only vertex with the same image as v, so $\varphi$ is actually
two to one on \six.

\subsection*{The quotient of \six}
\six\ has a fixed point free
involution given by sending {\it x} to {\it --x}, where we think of
\six\ as given by quaternions.  This quotient \three\ has 60 vertices,
300 cells, and is a triangulation of projective 3-space such that
the link of every vertex is an icosahedron.  We can consider \three\
as a higher dimensional analog of $K_5$ imbedded in the projective
plane, since that triangulation is the quotient of the icosahedron by
its involution.  The dual of \three\ is an analog of the Peterson graph.  

\subsection*{$B($\six$) = B($\three$)$}
The existence of the covering map $\six\longrightarrow\three$ 
shows that each coloring of \three\ determines a coloring of \six.
The vertices of $B(\six)$ are given by embedded 24-cells, and 
the 24-cells are fixed under the involution, so it follows that every
vertex of $B(\six)$ determines a unique vertex of $B(\three)$.  This
gives the desired isomorphism.

\subsection*{The map \three$\longrightarrow$\btwosix\ is 1-1.}
We saw above that the only distinct points $v,w$ satisfying
$\varphi(v) = \varphi(w)$ were $v=-w$.  This shows that the map from
\three\ to \btwosix\ is 1-1.  In other words,
\begin{quote} 
	There is a set of 300 5 by 5 Latin Squares whose structure is
isomorphic to \three.
\end{quote}

\subsection*{Labeling vertices of \six\ with permutations} 
The vertices of $S_5$ can be identified with the permutations on
$\{1,2,3,4,5\}$.  Using the map \six$\longrightarrow$\btwosix $= S_5$,
the vertices of \six\ can also be so identified.  Computationally, we
found that the 120 vertices determined 60 permutations, all the even
ones.  $v$ and $-v$ determine the same permutation.


\subsection*{\btwosix\ contains 4 disjoint copies of \three.}
Since the vertices of \btwosix\ are permutations, there is a
map given by $\sigma \longrightarrow \sigma^{-1}$ that is an
automorphism $\eta$ of the simplicial structure.  The two embeddings
of \three, $\varphi($\three$)$ and $\eta(\varphi($\three$))$, have no
4-simplices in common. If $\tau$ is any odd permutation, then $\tau
\cdot \varphi(p)$ and $\tau \cdot \eta(\varphi(p))$ give two more
copies of \three, whose vertex set is the set of odd permutations.

\subsection*{The structure of the remainder} 
	Once we have removed the 1200 4-simplices from the 1344
4-simplices of $S_5$ that are in the four copies of \three, we are
left with 144 4-simplices.  These also have an amazing structure.
They break into two isomorphic complexes with 60 vertices and 72
4-simplices, each isomorphic to the space of even permutations
$Alt_5$.  The vertices of this complex are all even permuations of
$1,2,3,4,5$.  Five permutations form a 4-simplex if it is possible to
make a Latin Square with them so that all rows and columns are even
permutations.

\subsection*{B$^2$(\six) is regular} 
Not only is it regular, but the link of a point is quite nice.  If
\igroup\ is the icosahedron, 
$$ Link(S_5,p) \approx B^2(\text{\igroup}). $$

\end{document}